
\input graphicx.tex
\input xy
\xyoption{all}
\magnification\magstephalf

\voffset0truecm
\hoffset=0truecm
\vsize=23truecm
\hsize=15.8truecm
\topskip=1truecm

\binoppenalty=10000
\relpenalty=10000

\font\tenbb=msbm10		\font\sevenbb=msbm7		\font\fivebb=msbm5
\font\tensc=cmcsc10		\font\sevensc=cmcsc7 	\font\fivesc=cmcsc5
\font\tensf=cmss10		\font\sevensf=cmss7		\font\fivesf=cmss5
\font\tenfr=eufm10		\font\sevenfr=eufm7		\font\fivefr=eufm5


\newfam\bbfam	\newfam\scfam	\newfam\frfam	\newfam\sffam

\textfont\bbfam=\tenbb
\scriptfont\bbfam=\sevenbb
\scriptscriptfont\bbfam=\fivebb

\textfont\scfam=\tensc
\scriptfont\scfam=\sevensc
\scriptscriptfont\scfam=\fivesc

\textfont\frfam=\tenfr
\scriptfont\frfam=\sevenfr
\scriptscriptfont\frfam=\fivefr

\textfont\sffam=\tensf
\scriptfont\sffam=\sevensf
\scriptscriptfont\sffam=\fivesf

\def\bb{\fam\bbfam \tenbb} 
\def\sc{\fam\scfam \tensc} 


\font\sezfont=cmbx10 scaled \magstep1
\font\subsectfont=cmbx10 scaled \magstephalf
\font\titfont=cmbx10 scaled \magstep2
\font\autfont=cmcsc10
\font\intfont=cmss10 

\let\no=\noindent
\let\bi=\bigskip
\let\me=\medskip
\let\sm=\smallskip
\let\ce=\centerline

\let\io=\infty
\def\qqquad{\quad\qquad}


\newcount\sectno\sectno=0
\newcount\subsectno\subsectno=0
\newcount\thmno\thmno=0
\newcount\tagno\tagno=0
\newcount\notitolo\notitolo=0
\newcount\defno\defno=0

\def\sect#1\par{
	\global\advance\sectno by 1 \global\subsectno=0\global\defno=0\global\thmno=0
	\vbox{\vskip.75truecm\advance\hsize by 1mm
	\hbox{\centerline{\sezfont \the\sectno.~~#1}}
	\vskip.25truecm}\nobreak}

\def\subsect#1\par{
	\global\advance\subsectno by 1
	\vbox{\vskip.75truecm\advance\hsize by 1mm
	\line{\subsectfont \the\sectno.\the\subsectno~~#1\hfill}
	\vskip.25truecm}\nobreak}
	
\def\defin#1{\global\advance\defno by 1
	\global\expandafter\edef\csname+#1\endcsname%
    {\number\sectno.\number\defno}
    \no{\bf Definition~\the\sectno.\the\defno.}}

\def\thm#1#2{
	\global\advance\thmno by 1
	\global\expandafter\edef\csname+#1\endcsname%
	{\number\sectno.\number\thmno}
	\no{\bf #2~\the\sectno.\the\thmno.}}

\def\Tag#1{\global\advance\tagno by 1 {(\the\tagno)}
    \global\expandafter\edef\csname+#1\endcsname%
    		{(\number\tagno)}}
\def\tag#1{\leqno\Tag{#1}}

\def\rf#1{\csname+#1\endcsname\relax}

\def\proof{\no{\sl Proof.}\enskip}
\def\qedn{\thinspace\null\nobreak\hfill\hbox{\vbox{\kern-.2pt\hrule height.2pt
        depth.2pt\kern-.2pt\kern-.2pt \hbox to2.5mm{\kern-.2pt\vrule
        width.4pt \kern-.2pt\raise2.5mm\vbox to.2pt{}\lower0pt\vtop
        to.2pt{}\hfil\kern-.2pt \vrule
        width.4pt\kern-.2pt}\kern-.2pt\kern-.2pt\hrule height.2pt
        depth.2pt \kern-.2pt}}\par\medbreak}
    \def\qed{\hfill\qedn}
    
\newif\ifpage\pagefalse
\newif\ifcen\centrue

\headline={
\ifcen\hfil\else
\ifodd\pageno
\global\hoffset=0.5truecm
\else
\global\hoffset=-0.4truecm
\fi\hfil
\fi}

\footline={
	\ifpage
		\hfill\rm\folio\hfill
	\else
		\global\pagetrue\hfill
\fi}

\lccode`\'=`\'

\def\bib#1{\me\item{[#1]\enskip}}

\def\ca#1{{\cal #1}}
\def\C{{\bb C}} \def\d{{\rm d}}
\def\R{{\bb R}} \def\Q{{\bb Q}}
 \def\N{{\bb N}}

\let\eps=\varepsilon 
\let\de=\partial

\mathchardef\void="083F

\def\diag{\mathop{\rm Diag}\nolimits}
\def\fix{\mathop{\rm Fix}\nolimits}

\def\id{\mathop{\rm Id}\nolimits}

\def\invlim{\mathop{\vtop{\offinterlineskip
\hbox{\rm lim}\kern1pt\hbox{\kern-1.5pt$\longleftarrow$}\kern-3pt}
}\limits}
\def\forevery#1#2{$$\displaylines{\hfilneg\rlap{$\qqquad\,\forall#1$}\hfil#2\cr}$$} 
\def\casi#1{\vcenter{\normalbaselines\mathsurround=0pt
		\ialign{$##\hfil$&\quad##\hfil\crcr#1\crcr}}}


\ce{\titfont Linearization of holomorphic germs} 

\ce{\titfont with quasi-Brjuno fixed points}

\me\ce{\autfont Jasmin Raissy}
\sm\ce{\intfont Dipartimento di Matematica, Universit\`a di Pisa}

\ce{\intfont Largo Bruno Pontecorvo 5, 56127 Pisa}

\sm\ce{\intfont E-mail: {\tt raissy@mail.dm.unipi.it}}
\bi

{\narrower

{\sc Abstract.} Let~$f$ be a germ of holomorphic diffeomorphism of~$\C^n$ fixing the origin~$O$, with~$\d f_O$ diagonalizable. We prove that, under certain arithmetic conditions on the eigenvalues of~$\d f_O$ and some restrictions on the resonances,~$f$ is locally holomorphically linearizable if and only if there exists a particular~$f$-invariant complex manifold. Most of the classical linearization results can be obtained as corollaries of our result. 

}
\bi

\sect Introduction 

We consider a germ of holomorphic diffeomorphism~$f$ of~$\C^n$ at a fixed point~$p$, which we may place at the origin~$O$. One of the main questions in the study of local holomorphic dynamics (see [A] and [B] for general surveys on this topic) is when~$f$ is {\it holomorphically linearizable}, i.e., when there exists a local holomorphic change of coordinates such that~$f$ is conjugated to its linear part. The answer to this question depends on the set of eigenvalues of~$\d f_O$, usually called the {\it spectrum} of~$\d f_O$. In fact if we denote by~$\lambda_1, \dots, \lambda_n\in \C^*$ the eigenvalues of~$\d f_O$, then it may happen that there exists a multi-index~$k=(k_1, \dots, k_n)\in \N^n$ with~$|k|=k_1+\cdots+k_n\ge 2$ and such that
$$\lambda^k - \lambda_j=\lambda_1^{k_1}\cdots\lambda_n^{k_n} - \lambda_j = 0\tag{eqres}$$
for some~$1\le j\le n$; a relation of this kind is called a {\it resonance} of~$f$, and~$k$ is called a {\it resonant multi-index}. A {\it resonant monomial} is a monomial~$z^k= z_1^{k_1}\cdots z_n^{k_n}$ in the~$j$-th coordinate, such that~$\lambda^k = \lambda_j$. From the formal point of view, we have the following classical result (see [Ar] pp.~192--193 for a proof):

\sm\thm{Te0.1}{Theorem} {\sl Let~$f$ be a germ of holomorphic diffeomorphism of~$\C^n$ fixing the origin~$O$ with no resonances. Then~$f$ is formally conjugated to its differential~$\d f_O$.}

\sm In presence of resonances, even the formal classification is not easy, as the following result of Poincar\'e-Dulac, [P], [D], shows

\sm\thm{Te0.2}{Theorem}(Poincar\'e-Dulac) {\sl Let~$f$ be a germ of holomorphic diffeomorphism of~$\C^n$ fixing the origin~$O$. Then~$f$ is formally conjugated to a formal power series~$g\in\C[\![z_1, \dots, z_n]\!]^n$ without constant term such that~$\d g_O$ is in Jordan normal form, and~$g$ has only resonant monomials.}

\sm The formal series~$g$ is called a {\it Poincar\'e-Dulac normal form} of~$f$; a proof of Theorem \rf{Te0.2} can be found in [Ar] p.~194. 

\sm Even without resonances, the holomorphic linearization is not guaranteed.  We need the following definitions:{\parindent=30pt
\sm\item{-} if all the eigenvalues of~$\d f_O$ have modulus less than~$1$, we say that the fixed point~$O$ is {\it attracting};
\sm\item{-} if all the eigenvalues of~$\d f_O$ have modulus greater than~$1$, we say that the fixed point~$O$ is {\it repelling};
\sm\item{-} if all the eigenvalues of~$\d f_O$ have modulus different from~$1$, we say that the fixed point~$O$ is {\it hyperbolic};
\sm\item{-} if all the eigenvalues of~$\d f_O$ are roots of unity, we say that the fixed point~$O$ is {\it parabolic}; in particular, if~$\d f_O = \id$ we say that~$f$ is {\it tangent to the identity};
\sm\item{-} if all the eigenvalues of~$\d f_O$ have modulus~$1$ but none is a root of unity, we say that the fixed point~$O$ is {\it elliptic};
\sm\item{-} if~$\d f_O=O$, we say that the fixed point~$O$ is {\it superattracting}.
\sm} 
The easiest positive result is due to Poincar\'e [P] who, using majorant series, proved the following

\sm\thm{Te0.3}{Theorem}(Poincar\'e, 1893 [P]) {\sl Let~$f$ be a germ of holomorphic diffeomorphism of~$\C^n$ with an attracting or repelling fixed point. Then~$f$ is holomorphically linearizable if and only if it is formally linearizable. In particular, if there are no resonances then~$f$ is holomorphically linearizable.}

\sm When~$O$ is not attracting or repelling, even without resonances, the formal linearization might diverge. Let us introduce the following definition:

\sm\defin{De1.2} Let~$n\ge2$ and let~$\lambda_1, \dots, \lambda_n\in\C^*$ be not necessarily distinct. Fix~$1\le s\le n$ and let~$\underline\lambda = (\lambda_{1},\dots, \lambda_s)$. For any~$m\ge 2$ put 
$$\omega_s(m) = \min_{2\le |k|\le m} \min_{1\le j\le n} |\underline \lambda^k - \lambda_j|,$$
where~$\underline \lambda^k= \lambda_{1}^{k_{1}}\cdots \lambda_s^{k_s}$. We say that~$\lambda=(\lambda_1, \dots, \lambda_n)$ {\it satisfies the partial Brjuno condition of order~$s$} if there exists a strictly increasing sequence of integers~$\{p_\nu\}_{\nu_\ge 0}$ with~$p_0=1$ such that
$$\sum_{\nu\ge 0} p_\nu^{-1} \log\omega_s(p_{\nu+1})^{-1}< \io.$$

\sm\thm{Re0.1}{Remark} For~$s=n$ the partial Brjuno condition of order~$s$ is nothing but the usual Brjuno condition introduced in [Br] (see also [M] pp.~25--37 for the one-dimensional case). When~$s<n$, the partial Brjuno condition of order~$s$ is indeed weaker than the Brjuno condition. Let us consider for example~$n=2$ and let~$\lambda, \mu\in\C^*$ be distinct. To check whether the pair~$(\lambda, \mu)$ satisfies the partial Brjuno condition of order~$1$, we have to consider only the terms~$|\lambda^k - \lambda|$ and~$|\lambda^k - \mu|$ for~$k\ge2$, whereas to check the full Brjuno condition we have to consider also the terms~$|\mu^h - \lambda|$,~$|\mu^h - \mu|$ for~$h\ge2$, and~$|\lambda^k\mu^h-\lambda|$,~$|\lambda^k\mu^h-\mu|$ for~$k,h\ge 1$. 

\sm\thm{Re0.2}{Remark} A~$n$-tuple~${\lambda= (\lambda_1, \dots, \lambda_s, 1,\dots, 1)\in(\C^*)^n}$ satisfies the partial Brjuno condition of order~$s$ if and only if~$(\lambda_1,\dots, \lambda_s)$ satisfies the Brjuno condition. 

\me We assume that the differential~$\d f_O$ is diagonalizable. Then, possibly after a linear change of coordinates, we can write
$$f(z) = \Lambda z + \hat f(z),$$
where~$\Lambda = \diag(\lambda_1, \dots,  \lambda_n)$, and~$\hat f$ vanishes up to first order at~$O\in \C^n$. 

The linear map~$z\mapsto\Lambda z$ has a very simple structure. For instance, for any subset~$\lambda_1, \dots, \lambda_s$ of eigenvalues with~$1\le s\le n$, the direct sum of the corresponding eigenspaces obviously is an invariant manifold on which this map acts linearly with these eigenvalues.  

\sm We have the following result of P\"oschel [P\"o] that generalizes the one of Brjuno [Br]:

\sm\thm{Te0.4}{Theorem}(P\"oschel, 1986 [P\"o]) {\sl Let~$f$ be a germ of holomorphic diffeomorphism of~$\C^n$ fixing the origin~$O$. If there exists a positive integer~$1\le s\le n$ such that the eigenvalues~$\lambda_1, \dots, \lambda_n$ of~$\d f_O$ satisfy the partial Brjuno condition of order~$s$, then there exists locally a complex analytic~$f$-invariant manifold~$M$ of dimension~$s$, tangent to the eigenspace of~$\lambda_{1},\dots,\lambda_s$ at the origin, on which the mapping is holomorphically linearizable.}

\sm In this paper we would like to extend P\"oschel Theorem in such a way to get a complete linearization in a neighbourhood of the origin.

\sm Before stating our result we need the following definitions:

\sm\defin{De1.1pre} Let~$1\le s\le n$. We say that~$\lambda = (\lambda_1, \dots, \lambda_s, \mu_1, \dots, \mu_r) \in(\C^*)^n$ {\it has only level~$s$ resonances} if there are only two kinds of resonances:
$$\lambda^k = \lambda_h \iff k\in\tilde K_1,$$
where
$$\tilde K_1 =\left \{k\in \N^n : |k|\ge 2, \sum_{p=1}^s k_p = 1~~\hbox{and}~~\mu_1^{k_{s+1}}\cdots \mu_r^{k_n}=1 \right\},$$
and
$$\lambda^k = \mu_j \iff k\in\tilde K_2,$$
where
$$\tilde K_2 = \{k\in \N^n : |k|\ge2, k_1=\cdots=k_s=0~\hbox{and}~\exists j \in\{1, \dots, r\}~\hbox{s.t.}~ \mu_1^{k_{s+1}}\cdots \mu_r^{k_n}=\mu_j \}.$$

\sm\thm{Ex1.3}{Example} When~$s<n$, if~$\lambda=(\lambda_1, \dots,\lambda_s, 1, \dots, 1)$ satisfies the Brjuno condition of order~$s$ then it is easy to verify that it has only level~$s$ resonances.

\sm\thm{Os1.1}{Remark} It is obvious that if the set~$\tilde K_2$ is empty (which implies that the set~$\tilde K_1$ is empty as well), there are no resonances. If~$\tilde K_1\ne\void$, having only level~$s$ resonances implies that the sets~$\{\lambda_1, \dots, \lambda_s\}$ and~$\{\mu_1, \dots, \mu_r\}$ are disjoint. If~$\tilde K_1=\void$ but~$\tilde K_2\ne\void$, then the sets~$\{\lambda_1, \dots, \lambda_s\}$ and~$\{\mu_1, \dots, \mu_r\}$ may intersect only in elements not involved in resonances, i.e., we can have~$\lambda_p=\mu_q$ for some~$p$ and~$q$ only if for every multi-index~$(k_{s+1},\dots, k_n)$, we have~$\mu_1^{k_{s+1}}\cdots\mu_r^{k_n}\ne\mu_q$, and for any resonance~$\mu_1^{k_{s+1}}\cdots\mu_r^{k_n}=\mu_j$ with~$j\ne q$, we have~$k_{s+q} = 0$.

\sm\thm{Ex1.1}{Example} Let~$\gamma\ge1$ and let~$\mu_3$ be a~$(\gamma+1)$-th primitive root of unity. Let~$\mu_1, \mu_2$ be two complex numbers of modulus different from~$1$ and such that
$$\mu_1^\alpha\mu_2^\beta = \mu_3$$
with~$\alpha,\beta\in\N\setminus\{0\}$. Then we have
$$\mu_1^\alpha\mu_2^\beta\mu_3^\gamma =1.$$ 
We can choose~$\mu_1, \mu_2$ such that the only resonant multi-indices for the triple~$(\mu_1, \mu_2, \mu_3)$ are~$(\alpha,\beta, 0)$,~$(\alpha-1,\beta, \gamma)$ and~$(\alpha,\beta-1,\gamma)$. Then, if we consider~$\lambda$ such that~$(\lambda, \mu_1, \mu_2, \mu_3)$ has only level~$1$ resonances, the admitted resonances are the following:
$$\eqalign{&\tilde K_1 = \{(1, \alpha,\beta,\gamma)\},\cr
			&\tilde K_2 = \{(0, \alpha,\beta,0), (0, \alpha-1,\beta, \gamma), (0, \alpha, \beta-1, \gamma)\}.}$$
			
\sm\thm{Ex1.2}{Example} Let us consider~$(\mu_1, \mu_2, \mu_3, \mu_4)\in(\C^*)^4$ with only one resonance, say~$\mu_1^p\mu_2^q=\mu_3$ with~$p,q\ge 1$, and such that~$(\lambda, \mu_1, \mu_2,\mu_3,\mu_4)$ has only level~$1$ resonances with~$\lambda=\mu_4$. Then
$$\eqalign{&\tilde K_1 = \void,\cr
			&\tilde K_2 = \{(0, p,q,0,0)\}.}$$
			
\sm\defin{De1.0} Let~$n\ge2$ and let~$\lambda_1, \dots, \lambda_n\in\C^*$ be not necessarily distinct. For any~$m\ge 2$ put 
$$\tilde\omega(m) = \min_{2\le |k|\le m\atop k\not\in {\ca Res_j(\lambda)}} \min_{1\le j\le n} |\lambda^k - \lambda_j|,$$
where~$\ca Res_j(\lambda)$ is the set of multi-indices~$k\in\N^n$ giving a resonance relation for~$\lambda =(\lambda_1, \dots, \lambda_n)$ relative to~$1\le j\le n$, i.e.,~$\lambda^k-\lambda_j=0$.
We say that~$\lambda$ {\it satisfies the reduced Brjuno condition} if there exists a strictly increasing sequence of integers~$\{p_\nu\}_{\nu_\ge 0}$ with~$p_0=1$ such that
$$\sum_{\nu\ge 0} p_\nu^{-1} \log\tilde\omega(p_{\nu+1})^{-1}< \io.$$

\sm\defin{De1.1} Let~$f$ be a germ of holomorphic diffeomorphism of~$\C^n$ fixing the origin~$O$ and let~$s\in\N$,~$1\le s\le n$. The origin~$O$ is called a {\it quasi-Brjuno fixed point of order~$s$} if~$\d f_O$ is diagonalizable and, denoting by~$\lambda=(\lambda_1, \dots, \lambda_n)$ the spectrum of~$\d f_O$, we have: {\parindent=30pt
\sm\item{(i)} $\lambda$ has only level~$s$ resonances;
\sm\item{(ii)} $\lambda$ satisfies the reduced Brjuno condition.
\sm}
\no We say that the origin is a {\it quasi-Brjuno fixed point} if there exists~$1\le s\le n$ such that it is a quasi-Brjuno fixed point of order~$s$.

\sm\defin{De0.1} Let~$f$ be a germ of holomorphic diffeomorphism of~$\C^n$ fixing the origin~$O$, and let~$1\le s\le n$. We will say that~$f$ {\it admits an osculating manifold~$M$ of codimension~$s$} if there is a germ of~$f$-invariant complex manifold~$M$ at~$O$ of codimension~$s$ such that the normal bundle~$N_M$ of~$M$ admits a holomorphic flat~$(1,0)$-connection that commutes with~$\d f|_{N_M}$.

\me We can now state our result which is a linearization result in presence of resonances:

\sm\thm{Teorema}{Theorem} {\sl Let~$f$ be a germ of a holomorphic diffeomorphism of~$\C^n$ having the origin~$O$ as a quasi-Brjuno fixed point of order~$s$. Then~$f$ is holomorphically linearizable if and only if it admits an osculating manifold~$M$ of codimension~$s$ such that~$f|_M$ is holomorphically linearizable.}

\me Roughly speaking, having only level $s$ resonances and the existence of the osculating manifold on which~$f$ is holomorphically linearizable take cares of the resonances in the $\mu_j$'s and give the formal linearization. Under these hypotheses the partial Brjuno condition of order~$s$ holds, so we have a partial holomorphic linearization given by P\"oschel's result, and the reduced Brjuno condition glues the formal linearization and the partial holomorphic linearization so to get a global holomorphic linearization.
In [R\"u], R\"ussmann gives an alternative way to pass from a formal linearization to a holomorphic one under an arithmetic hypothesis on the eigenvalues which implies the reduced Brjuno condition.

\me The structure of this paper is as follows. 

In the next section we shall explain the relations between the quasi-Brjuno condition and the partial Brjuno condition of order~$s$.

In the third section we shall give a characterization of osculating manifolds.

In the fourth section we shall prove a formal linearization result. 

In the fifth section we shall prove the holomorphic linearization result, i.e., Theorem \rf{Teorema}.
 
In the last section we shall point out similarities and differences with the known results.

\sm In the rest of the paper we shall denote by~$\|\cdot\|$ the norm~$\|\cdot\|_\io$; but we could also had used the norm~$\|\cdot\|_2$ thanks to the equivalence of such norms. We shall also need the following notation: if~$g\colon \C^n\to \C$ is a holomorphic function with~$g(O)=0$ (or a formal power series without constant term), and~$z=(x, y)\in \C^n$ with~$x\in\C^s$ and~$y\in\C^{n-s}$, we shall denote by~${\rm ord}_x(g)$ the maximum positive integer~$m$ such that~$g$ belongs to the ideal~$(x_1,\cdots, x_s)^m$.

\sect Quasi-Brjuno condition vs Partial Brjuno condition

Notice that whereas it is always possible to introduce the reduced Brjuno condition, the partial Brjuno condition makes sense only when there are no resonant multi-indices~$k\in\N^n$, with~$|k|\ge 2$ and~$k_{s+1}=\dots=k_n=0$.
Anyway, when we have only level~$s$ resonance, we can deal with these two condition at the same time.

\sm\thm{Re1.1pre}{Remark} If~$\lambda$ has only level~$s$ resonances, then we have
$$\tilde\omega(m)= \min_{2\le |k|\le m}\min\left\{\min_{1\le j\le n\atop k_1+\cdots+k_s\ge2}|\lambda^k-\lambda_j|, \min_{1\le j\le n-s\atop k_1+\cdots+k_s=1}|\lambda^k-\lambda_{s+j}|\right\},$$
therefore
$$\tilde\omega(m)= \min\left\{\omega_s(m), \min_{2\le |k|\le m\atop (k_{s+1},\dots, k_n)\ne O}\left\{\min_{1\le j\le n\atop k_1+\cdots+k_s\ge2}|\lambda^k-\lambda_j|, \min_{1\le j\le n-s\atop k_1+\cdots+k_s=1}|\lambda^k-\lambda_{s+j}|\right\}\right\},$$
so it is obvious that, since~$\tilde\omega(m)\le \omega_s(m)$ for every~$m\ge 2$, the reduced Brjuno condition implies the partial Brjuno condition of order~$s$. A partial converse is the following

\sm\thm{Le0.1}{Lemma} {\sl Let~$n\ge2$ and let~$\lambda_1, \dots, \lambda_n\in\C^*$ be not necessarily distinct. Let~$1\le s\le n$ be such that~$\lambda=(\lambda_1, \dots, \lambda_n)$ has only level~$s$ resonances. Then, if there exists a strictly increasing sequence of integers~$\{p_\nu\}_{\nu_\ge 0}$ with~$p_0=1$ such that
$$\sum_{\nu\ge 0} p_\nu^{-1} \log\omega_s(p_{\nu+1})^{-1}< \io,$$
(i.e.,~$\lambda$ satisfies the partial Brjuno condition of order~$s$), and there exist~$k\in \N$ and~$\alpha\ge 1$ such that
$$p_\nu> k\Rightarrow {\tilde \omega(p_\nu-k) \ge \omega_s(p_\nu)^\alpha,}$$
then~$\lambda$ satisfies the reduced Brjuno condition.}

\sm\proof Let~$q_0=p_0$ and~$q_j=p_{\nu_0+j} - k$ for~$j\ge 1$, where~$\nu_0$ is the minimum index such that~$p_\nu>k$ for all~$\nu\ge \nu_0$. Then we have
$$\eqalign{\sum_{\nu\ge 0} q_\nu^{-1} \log\tilde\omega(q_{\nu+1})^{-1} &\le \alpha\sum_{\nu\ge 0} q_\nu^{-1}\log\omega_s(q_{\nu+1}+k)^{-1}\cr
	&= \alpha p_0^{-1}\log\omega_s(p_{\nu_0+1})^{-1} + \alpha\sum_{\nu\ge \nu_0+2}{p_\nu\over p_\nu -k} p_\nu^{-1} \log\omega_s(p_{\nu+1})^{-1}\cr
	&\le 2\alpha\sum_{\nu\ge 0} p_\nu^{-1} \log\omega_s(p_{\nu+1})^{-1}\cr
	&< \io,}$$
and we are done. \qed

\sm\thm{Re1.1}{Remark} Suppose that~$\lambda$ has only level~$s$ resonances. Recall that a sequence~$\{a_m\}$ is said to be {\it Diophantine of exponent~$\tau>1$} if there exist~$\gamma$,~$\gamma'> 0$ so that~$\gamma' m^{-\beta}\ge a_m\ge\gamma m^{-\beta}$ (see also [C], [G] and [S]).
Then if~$\tilde\omega(m)$ is Diophantine of exponent~$\beta>1$, and if~$\omega_s(m)$ is Diophantine of exponent~$\eps>1$, there always exist~$\alpha\ge 1$ and~$\delta>0$ for which
$$\tilde\omega(m)\ge\gamma m^{-\beta}\ge \delta m^{-\eps\alpha} \ge \omega_s(m)^\alpha,$$
and thus the hypothesis of Lemma \rf{Le0.1} is satisfied with~$k=0$.

More in general, if we have
\forevery{m\ge k+2}{\tilde\omega(m-k)\ge \omega_s(m)^\alpha}
for some~$k\in \N$ and~$\alpha\ge 1$, the hypothesis of Lemma \rf{Le0.1} is obviously satisfied. 
For example if~$\lambda_1,\dots, \lambda_s\in\R$ are positive and~$\lambda_{s+1}, \dots,\lambda_n\in \{-1, +1\}$ then it is easy to verify that
\forevery{m\ge 3}{\tilde\omega(m-1)\ge\omega_s(m).}
Furthermore, if~$\lambda_{s+1}= \cdots=\lambda_n=1$ then~$\tilde\omega(m)=\omega_s(m)$, and so in this case the partial Brjuno condition of order~$s$ coincides with the reduced Brjuno condition.

\sect Osculating manifolds

Let~$f$ be a germ of holomorphic diffeomorphism of~$\C^n$ at a point which we may assume without loss of generality to be the origin~$O$, and let~$M$ be an~$f$-invariant complex manifold through~$O$ of codimension~$s$, with~$1\le s\le n$. In this situation, the differential~$\d f$ acts on the normal bundle~$N_M= T\C^n/TM$. 

\sm It is obvious that locally every holomorphic bundle admits a holomorphic flat~$(1,0)$-connection (it suffices to take the trivial connection on a trivialization).  Moreover, it is easy to prove the following result, which has exactly the same proof as in the smooth case (adopting for instance the argument in [BCS] pp. 272--274).

\sm\thm{Prcoordinate}{Proposition} {\sl Let~$\pi\colon E\to M$ be a holomorphic vector bundle on a complex manifold~$M$ and let~$\nabla$ be a holomorphic flat~$(1,0)$-connection. Then there are a local holomorphic coordinate system about~$O$ and a local holomorphic frame of~$E$ in which all the connection coefficients~$\Gamma^i_{jk}$ are zero.}

\sm In the particular case of the normal bundle we have the following useful result.

\sm\thm{LeNormale}{Lemma} {\sl Let~$M\subset \C^n$ be a complex manifold through~$O$ of codimension~$s$, with~$1\le s\le n$ and let~$N_M$ be its normal bundle. Fix~$p\in M$. Take a local holomorphic frame in a neighbourhood of~$p$. Then there exist local coordinates at~$p$ in~$\C^n$ such that for every local holomorphic frame~$\{V_1, \dots, V_s\}$ of~$N_M$ we can find local holomorphic coordinates~$(U,Z)$ with~$z=(x,y)$, adapted to~$M$ (i.e.,~$M\cap U=\{x=0\}$) such that, on~$U\cap M$,
$$ V_j = \pi\left({\de \over\de x_j}\right)$$
for every~$j=1, \dots, s$, where~$\pi\colon T\C^n\to N_M$ is the canonical projection.}

\sm\proof Let us choose local holomorphic coordinates~$\tilde z = (\tilde x, \tilde y)$ at~$p$ adapted to~$M$. Then for every point~$(0, \tilde y)\in M$ there exists a non-singular matrix~$A(\tilde y) = (a_{ij}(\tilde y))$, depending holomorphically on~$\tilde y$, such that
$$V_j(\tilde y) = \sum_{i=1}^s a_{ij}(\tilde y) \pi\left.\left({\de\over \de \tilde x_i}\right)\right|_{(0, \tilde y)}.$$
Therefore, using the coordinates 
$$\eqalign{&x_i=\sum_{i=1}^s a_{ij}(\tilde y) \tilde x_i \quad\hbox{for}~i=1,\dots, s,\cr
  		&y_j=\tilde y_j \quad\quad\quad\hbox{for}~j=1,\dots, r,}$$
we obtain the assertion. \qed

\me\defin{Deconnessione} Let~$f$ be a germ of holomorphic diffeomorphism of~$\C^n$ fixing the origin~$O$, and let~$M$ be a germ of~$f$-invariant complex manifold at~$O$ of codimension~$s$, with~$1\le s\le n$. We say that a holomorphic flat~$(1,0)$-connection~$\nabla$ of the normal bundle~$N_M$ of~$M$ is {\it $f$-invariant} if it commutes with~$\d f|_{N_M}$.



\sm\thm{PrOsculating}{Theorem} {\sl Let~$f$ be a germ of holomorphic diffeomorphism of~$\C^n$ fixing the origin~$O$, let~$M$ be a germ of~$f$-invariant complex manifold through~$O$ of codimension~$s$, with~$1\le s\le n$, and let~$\nabla$  be a holomorphic flat~$(1,0)$-connection of the normal bundle~$N_M$. Then~$\nabla$ is~$f$-invariant if and only if there exist local holomorphic coordinates~$z=(x,y)$ about~$O$ adapted to~$M$ in which~$f$ has the form
$$\eqalign{&x_i'=\lambda_i x_i +\eps_i x_{i+1}+ f^1_i(x,y) \quad\hbox{for}~i=1,\dots, s,\cr
  		&y_j'=\mu_j y_j +\eps_{s+j} y_{j+1}+ f^2_j(x,y) \quad\hbox{for}~j=1,\dots, r=n-s,}\tag{obsc}$$
where~$\eps_i,\eps_{s+j}\in\{0,1\}$, and
$${\rm ord}_x(f_i^1)\ge 2,$$
for any~$i=1, \dots, s$.}
			
\sm\proof If there exist local holomorphic coordinates~$z=(x,y)$ about~$O$ adapted to~$M$, in which~$f$ has the form \rf{obsc} with~${\rm ord}_x(f_i^1)\ge 2$ for any~$i=1, \dots, s$, then it is obvious to verify that the trivial holomorphic flat~$(1,0)$-connection is~$f$-invariant.

Conversely, let~$\nabla$ be a holomorphic flat~$f$-invariant~$(1,0)$-connection of the normal bundle~$N_M$.
Thanks to Proposition \rf{Prcoordinate} and to Lemma \rf{LeNormale} we can find local holomorphic coordinates~$z=(x,y)$ adapted to~$M$, in which all the connection coefficients~$\Gamma^i_{jk}$ with respect to the local holomorphic frame~$\{\pi({\de\over\de x_1}), \dots, \pi({\de\over\de x_s})\}$ of~$N_M$ are zero.
We may assume without loss of generality, (up to linear changes of the coordinates we can assume that the linear part of~$f$ is in Jordan normal form), that in such coordinates~$f$ has the form
$$\eqalign{&x_i'=\lambda_i x_i+ \eps_i x_{i+1}+ f^1_i(x,y) \quad\hbox{for}~i=1,\dots, s,\cr
  		&y_j'=\mu_j y_j +\eps_{s+j} y_{j+1}+ f^2_j(x,y) \quad\hbox{for}~j=1,\dots, r,}$$
where~$\eps_i,\eps_{s+j}\in\{0,1\}$. Moreover, since~$M=\{x=0\}$ is~$f$-invariant, we have
$${\rm ord}_x(f_i^1)\ge 1.$$
Thanks to the~$f$-invariance of~$\nabla$ we have
$$\displaystyle{\nabla_{\de\over\de y_k} \left(\d f|_{N_M}\pi\left({\de\over \de x_j}\right)\right) = \d f|_{N_M}\nabla_{\de\over\de y_k} \pi\left({\de\over \de x_j}\right)}$$
for any~$j=1, \dots, s$ and~$k=1, \dots, r$. Now the right-hand side vanishes, because in the chosen coordinates we have~$\nabla_{\de\over\de y_k} \pi\left({\de\over \de x_j}\right)= 0$. So, using Leibniz formula, we obtain
$$\eqalign{0 &=\nabla_{\de\over\de y_k} \left(\d f\pi\left({\de\over \de x_j}\right)\right) \cr 
&=\nabla_{\de\over\de y_k} \left(\sum_{h=1}^s\left(\lambda_h\delta_{hj} + \eps_h\delta_{h,j+1} + {\de f^1_h\over \de x_j}(0,y)\right)\pi\left({\de\over\de x_h}\right)\right)\cr
&= \sum_{h=1}^s\left(\lambda_h\delta_{hj} + \eps_h\delta_{h,j+1} + {\de f^1_h\over \de x_j}(0,y)\right)\nabla_{\de\over\de y_k} \pi\left({\de\over \de x_h}\right) + \sum_{h=1}^s {\de \over\de y_k}\left({\de f^1_h \over\de x_j}(0,y)\right)\pi\left({\de \over \de x_h}\right)\cr
&=\sum_{h=1}^s {\de \over\de y_k}\left({\de f^1_h \over\de x_j}(0,y)\right)\pi\left({\de \over \de x_h}\right).}\tag{latosinistro}$$
Therefore we obtain
$${\de \over\de y_k}\left({\de f^1_h \over\de x_j}(0,y)\right)=0$$
for every~$j,h=1,\dots, s$ and~$k=1, \dots r$, that is
$${\rm ord}_x(f_h^1)\ge 2$$
for every~$h=1,\dots, s$, and this concludes the proof. \qed

\sm\thm{CoObs}{Corollary} {\sl Let~$f$ be a germ of holomorphic diffeomorphism of~$\C^n$ fixing the origin~$O$, and let~$1\le s\le n$. Then~$f$ admits an osculating manifold~$M$ of codimension~$s$ such that~$f|_M$ is holomorphically linearizable if and only if there exist local holomorphic coordinates~$z=(x,y)$ about~$O$ adapted to~$M$ in which~$f$ has the form
$$\eqalign{&x_i'=\lambda_i x_i +\eps_i x_{i+1}+ f^1_i(x,y) \quad\hbox{for}~i=1,\dots, s,\cr
  		&y_j'=\mu_j y_j +\eps_{s+j} y_{j+1}+ f^2_j(x,y) \quad\hbox{for}~j=1,\dots, r,}\tag{obsc3}$$
where~$\eps_i,\eps_{s+j}\in\{0,1\}$, and
$$\eqalign{&{\rm ord}_x (f_i^1)\ge 2,\cr
  			&{\rm ord}_x (f_j^2)\ge 1, }\tag{obsc2}$$
for any~$i=1, \dots, s$ and~$j=1,\dots, r$.}

\sm\proof One direction is clear. Conversely, thanks to Theorem \rf{PrOsculating}, the fact that~$M$ is osculating, i.e.,~$M$ is an~$f$-invariant complex manifold through~$O$ of codimension~$s$, with~$1\le s\le n$, with a holomorphic flat~$f$-invariant~$(1,0)$-connection of the normal bundle~$N_M$, is equivalent to the existence of local holomorphic coordinates~$z=(x,y)$ about~$O$ adapted to~$M$, in which~$f$ has the form \rf{obsc3}\ with~${\rm ord}_x (f_i^1)\ge 2$ for any~$i=1, \dots, s$. 

%
Furthermore,~$f|_M$ is linearizable; therefore there exists a local holomorphic change of coordinate, tangent to the identity, and of the form
$$\eqalign{&\tilde x= x\cr
  		   &\tilde y = \Phi(y),}$$
conjugating~$f$ to~$\tilde f$ of the form \rf{obsc3} satisfying \rf{obsc2}, as we wanted. \qed

Then we could say that, if we write~$f$ as in \rf{obsc}, the hypothesis of~$f$-invariance is equivalent to~${\rm ord}_x (f_i^1)\ge 1$; $f|_M$ linearized is equivalent to~${\rm ord}_x (f_j^2)\ge 1$; osculating means that~$f_i^1$ has no terms of order~$1$ in~$x$, that is,~$f^1_i = \sum_{h,k}x_h x_k \theta_i^{hk}(x,y)$.

Notice that in Theorem \rf{PrOsculating} and in Corollary \rf{CoObs}, up to linear changes of coordinates, we can always assume~$\eps_i,\eps_j\in\{0,\eps\}$ instead of~$\eps_i,\eps_j\in\{0,1\}$ for every~$\eps>0$ small enough.

\me Since we are going to first prove a formal result, we need the formal analogue of Definition~\rf{De0.1}. We define a {\it formal complex manifold~$M$ of codimension~$s$} by means of an ideal of formal complex power series generated by~$s$ power series~$g_1, \dots, g_s$ such that their differentials at the origin~$\d g_1, \dots \d g_s$ are linearly independent (see also [BER] and [BMR]). Denote by~$\widehat{T\C^n}$ the formal tangent bundle of~$\C^n$, that is the space of all formal vector fields with complex coefficients. Then the {\it formal tangent bundle~$\widehat{T M}$ to~$M$} is well-defined as being the set of formal vector fields of~$\widehat{T\C^n}$ vanishing on the ideal of formal power series generated by~$g_1, \dots, g_s$. The {\it formal normal bundle~$\widehat{N_M}$ of~$M$} is then the quotient~$\widehat{T\C^n}/\widehat{TM}$. A {\it formal connection} on the formal normal bundle is a formal map~$\widehat \nabla\colon \widehat{TM}\times \widehat{N_M} \to \widehat{N_M}$ which satisfies the usual properties of a connection but in the formal category. Thus the following definitions makes sense.

\me\defin{Deconnessioneformale} Let~$f$ be a formal invertible map of~$\C^n$ without constant term, and let~$M$ be an~$f$-invariant formal complex manifold of codimension~$s$, with~$1\le s\le n$. We say that a formal flat~$(1,0)$-connection~$\widehat\nabla$ of the formal normal bundle~$\widehat{N_M}$ of~$M$ is {\it $f$-invariant} if it commutes with~$\d f|_{N_M}$.

\sm\defin{De0.1formal} Let~$1\le s\le n$, and let~$f$ be a formal invertible map of~$\C^n$ without constant term. We will say that~$f$ {\it admits a formal osculating manifold~$M$ of codimension~$s$} if there is an~$f$-invariant formal complex manifold~$M$ of codimension~$s$ such that the formal normal bundle~$\widehat{N_M}$ of~$M$ admits a formal flat~$f$-invariant~$(1,0)$-connection.

\sm Then, for the formal normal bundle we can prove the formal analogue of Proposition \rf{Prcoordinate} (using a formal solution of the parallel transport equation that can be easily computed) and Lemma \rf{LeNormale}. We then have the following results, whose proofs are the formal analogues of the ones of Theorem \rf{PrOsculating} and Corollary \rf{CoObs}.

\sm\thm{PrOsculatingformal}{Theorem} {\sl Let~$f$ be a formal invertible map of~$\C^n$ without constant term, let~$M$ be an~$f$-invariant formal complex manifold through~$O$ of codimension~$s$, with~$1\le s\le n$, and let~$\widehat\nabla$ be a formal flat~$(1,0)$-connection of the formal normal bundle~$\widehat{N_M}$. Then~$\widehat\nabla$ is~$f$-invariant if and only if there exist local formal coordinates~$z=(x,y)$ about~$O$ adapted to~$M$ in which~$f$ has the form
$$\eqalign{&x_i'=\lambda_i x_i +\eps_i x_{i+1}+ f^1_i(x,y) \quad\hbox{for}~i=1,\dots, s,\cr
  		&y_j'=\mu_j y_j +\eps_{s+j} y_{j+1}+ f^2_j(x,y) \quad\hbox{for}~j=1,\dots, r,}\tag{obscfor}$$
where~$\eps_i,\eps_{s+j}\in\{0,1\}$, and
$${\rm ord}_x (f_i^1)\ge 2,$$
for any~$i=1, \dots, s$.}

\sm\thm{CoObsform}{Corollary} {\sl Let~$1\le s\le n$, and let~$f$ be a formal invertible map of~$\C^n$ without constant term. Then~$f$ admits a formal osculating manifold~$M$ of codimension~$s$ such that~$f|_M$ is formally linearizable if and only if there exist local formal coordinates~$z=(x,y)$ about~$O$ adapted to~$M$ in which~$f$ has the form
$$\eqalign{&x_i'=\lambda_i x_i +\eps_i x_{i+1}+ f^1_i(x,y) \quad\hbox{for}~i=1,\dots, s,\cr
  		&y_j'=\mu_j y_j +\eps_{s+j} y_{j+1}+ f^2_j(x,y) \quad\hbox{for}~j=1,\dots, r,}\tag{obsc3for}$$
where~$\eps_i,\eps_{s+j}\in\{0,1\}$, and
$$\eqalign{&{\rm ord}_x(f_i^1)\ge 2,\cr
  			&{\rm ord}_x(f_j^2)\ge 1, }\tag{obsc2for}$$
for any~$i=1, \dots, s$ and~$j=1,\dots, r$.}

\sect Formal linearization

As announced, we first prove a formal result.

\sm\thm{Te1.0}{Theorem} {\sl Let~$f$ be a formal invertible map of~$\C^n$ without constant term such that~$\d f_O$ is diagonalizable and the spectrum of~$\d f_O$ has only level~$s$ resonances, with~$1\le s\le n$. Then~$f$ is formally linearizable if and only if it admits an osculating formal manifold of codimension~$s$ such that~$f|_M$ is formally linearizable.}

\sm\proof If~$f$ is formally linearizable the assertion is obvious. 

Conversely, using Corollary \rf{CoObsform}, we can choose formal local coordinates
$$(x,y)=(x_1, \dots, x_s, y_1, \dots, y_r)$$
such that, writing~$(x',y') = f(x,y)$,~$f$ is of the form
$$\eqalign{&x_i'=\lambda_i x_i + f^1_i(x,y) \quad\hbox{for}~i=1,\dots, s,\cr
  			&y_j'=\mu_j y_j + f^2_j(x,y) \quad\hbox{for}~j=1,\dots, r,}$$
where
$$\eqalign{&{\rm ord}_x(f_i^1)\ge 2,\cr
  			&{\rm ord}_x(f_j^2)\ge 1. }$$
Denote by~$\Lambda$ the diagonal matrix~$\diag(\lambda_1,\dots, \lambda_s, \mu_1,\dots, \mu_r)$. We would like to prove that a formal solution~$\psi$ of
$$f\circ \psi=\psi\circ \Lambda\tag{eq1}$$
exists of the form
$$\eqalign{&x_i= u_i + \psi^1_i(u,v) \quad\hbox{for}~i=1,\dots, s,\cr
  			&y_j= v_j + \psi^2_j(u,v) \quad\hbox{for}~j=1,\dots, r, }$$
where~$(u,v)=(u_1, \dots, u_s,v_1, \dots,v_r)$ and~$\psi^1_i$ and~$\psi^2_j$ are formal power series with
$$\eqalign{&{\rm ord}_u(\psi_i^1)\ge 2,\cr
  			&{\rm ord}_u(\psi_j^2)\ge 1. }$$

Write~$f(z)=\Lambda z + \hat f(z)$ and~$\psi(w) = w + \hat \psi(w)$, where~$z=(x,y)$ and~$w=(u,v)$. Then equation \rf{eq1} is equivalent to
$$\hat\psi\circ \Lambda - \Lambda \hat\psi = \hat f\circ \psi. \tag{eq2}$$
To obtain a formal solution, we first write
$$\hat\psi = \sum_{|k|\ge 2} \psi_k w^k, \quad \psi_k \in \C^n,$$
where~$k=(k_1,\dots, k_n)$, and
$$\hat f = \sum_{|l|\ge 2} f_l z^l, \quad f_l \in \C^n,$$
where~$l=(l_1,\dots, l_n)$. Denoting~$\tilde\lambda=(\lambda_1, \dots, \lambda_s, \mu_1,\dots,\mu_r)=(\tilde\lambda_1, \dots, \tilde\lambda_n)$, equation \rf{eq2} becomes
$$\sum_{|k|\ge 2} A_k \psi_k w^k = \sum_{|l|\ge 2}f_l \left(\sum_{|m|\ge1}\psi_m w^m \right)^l, \tag{eq3}$$
where
$$A_k=\tilde\lambda^k I_n-\Lambda.$$
The matrices~$A_k$ might not be invertible for some choice of~$k$ due to the presence of resonances. We can write~$A_k=\diag(A^1_k, A^2_k)$ and recall that having only level~$s$ resonances means that~$\det(A_k^1) = 0$ if and only if
$$k\in\tilde K_1,$$
and~$\det(A_k^2) = 0$ if and only if
$$k\in\tilde K_2.$$
Moreover, from the hypotheses of the Theorem we have that~$f^1_l=0$ for~$l$ in~$K_1\cup K_2$ and~$f^2_l=0$ for~$l$ in~$K_2$, where
$$\eqalign{&K_1=\{l\in \N^n : |l|\ge2, l=( 0, \dots, 0,l_{i},0,\dots,0, l_{s+1}, \dots, l_n), l_{i}=1~\hbox{and}~i\in\{1, \dots, s\}\}\cr
& K_2=\{l\in \N^n : |l|\ge2, l=(0, \dots, 0, l_{s+1}, \dots, l_n)\}.}$$ 
Notice that~$\tilde K_1\subseteq K_1$ and~$\tilde K_2\subseteq K_2$. For each~$j$ in~$\{1, \dots, s\}$, let us denote by~$K_1^j$ the set $\{l\in \N^n : |l|\ge2, l=( 0, \dots, 0,l_{j},0,\dots,0, l_{s+1}, \dots, l_n), l_{j}=1\}$, so that~$K_1 = \cup_{j=1}^s K_1^j$. We look for a solution of \rf{eq1} with~$\psi^1_k=0$ for~$k\in K_1\cup K_2$ and~$\psi^2_k=0$ for~$k\in K_2$. 

To do so, let us write \rf{eq3} in a more explicit way: for~$i=1,\dots, s$
$$\sum_{|k|\ge 2\atop k\not\in K_1\cup K_2} (\tilde\lambda^k - \lambda_i) \psi^1_{k,i} w^k = \sum_{|l|\ge 2\atop l\not\in K_1\cup K_2}f^1_{l,i} \left(\sum_{|m|\ge1}\psi_m w^m \right)^l,\tag{eq3bbis}
$$
and for~$j=1,\dots, r$
$$\eqalign{\sum_{p=1}^s\sum_{|k|\ge 2\atop k\in K_1^p} (\tilde\lambda^k - \mu_j) \psi^2_{k,j}& w^k + \sum_{|k|\ge 2\atop k\not\in K_1\cup K_2} (\tilde\lambda^k - \mu_j) \psi^2_{k,j} w^k \cr 
&= \sum_{p=1}^s\sum_{|l|\ge 2\atop l\in K^p_1}f^2_{l,j} \left(\sum_{|m|\ge1}\psi_m w^m \right)^l + \sum_{|l|\ge 2\atop l\not\in K_1\cup K_2}f^2_{l,j} \left(\sum_{|m|\ge1}\psi_m w^m \right)^l.
}\tag{eq3bis}$$

Now, it is obvious that there are no terms~$w^k$ with~$k\in K_2$ in either side of \rf{eq3bbis} and of \rf{eq3bis}, and we can obtain terms~$w^k$ with~$k\in K_1$ in \rf{eq3bis} only from terms with~$l\in K_1$. In fact, if~$l\in K_1^h$ then
$$\eqalign{\left(\sum_{|m|\ge1}\psi_m w^m \right)^l &=\left(u_{h} + \sum_{p, q}u_p u_q \theta_h^{pq}(u,v) \right) \left(\prod_{j=1}^r \big(v_j + \sum_{p}u_p \theta_j^{p}(u,v) \big)^{l_{s+j}} \right)\cr
				&= u_h v_1^{l_{s+1}}\cdots v_r^{l_n} + \sum_{p, q}u_p u_q \chi^{pq}(u,v) \cr
				&= w^l + \sum_{p, q}u_p u_q \chi^{pq}(u,v).
}$$
Therefore for~$j=1,\dots, r$, we have 
$$\eqalign{\sum_{p=1}^s\sum_{|k|\ge 2\atop k\in K_1^p} (\tilde\lambda^k - \mu_j) \psi^2_{k,j} w^k &= \sum_{p=1}^s\sum_{|l|\ge 2\atop l\in K^p_1}f^2_{l,j} \left(\sum_{|m|\ge1}\psi_m w^m \right)^l \cr 
				&= \sum_{p=1}^s\sum_{|l|\ge 2\atop l\in K^p_1}f^2_{l,j}\big( w^l + \sum_{a, b}u_a u_b \chi^{ab}(u,v)\big)
				}$$
from which we conclude that for~$k\in K_1^p$ and~$j=1, \dots, r$ we have
$$\psi^2_{k,j} = {f^2_{k,j} (\tilde\lambda^k - \mu_j)^{-1}}.\tag{eqK_2}$$
The remaining~$\psi_k$ with~$k\not\in K_1\cup K_2$ are easily determined by recursion, as usual. \qed

\sect Holomorphic linearization

Now we can prove the main result of this paper.

\sm\thm{Te1.1}{Theorem} {\sl Let~$f$ be a germ of a holomorphic diffeomorphism of~$\C^n$ having the origin~$O$ as a quasi-Brjuno fixed point of order~$s$, with~$1\le s\le n$. Then~$f$ is holomorphically linearizable if and only if it admits an osculating manifold~$M$ of codimension~$s$ such that~$f|_M$ is holomorphically linearizable.}

\sm\proof If~$f$ is linearizable the assertion is obvious. 

Conversely, we already know, thanks to the previous result, that~$f$ is formally linearizable, (notice that, thanks to Corollary \rf{CoObs}, the changes of coordinates needed before finding~$\psi$ are holomorphic because now~$M$ is a complex manifold).
Since the spectrum of~$\d f_O$ satisfies the reduced Brjuno condition, to prove the convergence of the formal conjugation~$\psi$ in a neighbourhood of the origin it suffices to show that
$$\sup_k{1\over |k|} \log\|\psi_k\|<\io .\tag{eq6}$$
Since~$f$ is holomorphic in a neighbourhood of the origin, there exists a positive number~$P$ such that~$\|f_l\|\le P^{|l|}$ for~$|l|\ge 2$. The functional equation \rf{eq1} remains valid under the linear change of coordinates~$f(z)\mapsto sf(z/Q)$,~$\psi(w)\mapsto Q\psi(w/Q)$ with~$Q=\max\{1, P^2\}$. Hence we may assume that 
\forevery{|l|\ge 2}{\|f_l\|\le 1.}
It follows from \rf{eq3} and \rf{eqK_2} that
$$\|\psi_k\|\le \cases{\displaystyle \eps_k^{-1} \sum_{k_1+\cdots +k_\nu= k \atop \nu\ge 2} \|\psi_{k_1}\|\cdots \|\psi_{k_\nu}\|, &$|k|\ge 2, \quad k\not\in K_1\cup K_2$, \cr\noalign{\sm}
\displaystyle \eps_k^{-1}, &$|k|\ge 2, \quad k\in K_1$,} \tag{eq7}$$
where
$$\eps_k = \cases{\displaystyle\min_{1\le i \le n} |\tilde\lambda^k - \tilde \lambda_i|,&$k\not\in K_1\cup K_2$,\cr\noalign{\sm}
\displaystyle\min_{1\le h \le r}|\tilde\lambda^k -\mu_h|, &$k\in K_1$.}$$
We can define, inductively, for~$j \ge 2$ 
$$\alpha_j= \sum_{j_1+\cdots + j_\nu =j \atop \nu \ge 2} \alpha_{j_1} \cdots \alpha_{j_\nu},$$
and for~$|k|\ge 2$
$$\delta_k = \cases{\displaystyle \eps_k^{-1}\max_{k_1+\cdots + k_\nu =k\atop \nu\ge 2} \delta_{k_1}\cdots\delta_{k_\nu},&$k\not\in K_1\cup K_2$,\cr\noalign{\sm}\displaystyle \eps_k^{-1}, &$k\in K_1$, \cr\noalign{\sm}
\displaystyle 0,&$k\in K_2$,}$$
with~$\alpha_1 =1$ and~$\delta_e= 1$, where~$e$ is any integer vector with~$|e|=1$. Then, by induction, we have that
\forevery{|k|\ge 1}{\|\psi_k\|\le \alpha_{|k|}\delta_k.}
Therefore, to establish \rf{eq6}, it suffices to prove analogous estimates for~$\alpha_j$ and~$\delta_k$.

\sm It is easy to estimate~$\alpha_j$. Let~$\alpha= \sum_{j \ge 1}\alpha_j t^j$. We have
$$\eqalign{\alpha - t &= \sum_{j\ge 2} \alpha_j t^j \cr
					&= \sum_{j\ge 2}\left(\sum_{h\ge 1}\alpha_h t^h\right)^j\cr
					&= {\alpha^2 \over 1- \alpha}.}$$
This equation has a unique holomorphic solution vanishing at zero
$$\alpha= {t+1 \over 4} \left(1 - \sqrt{1-{8t\over(1+t)^2}}\right),$$
defined for~$|t|$ small enough. Hence,
$$\sup_j {1\over j}\log \alpha_j < \io,$$
as we want.

\sm To estimate~$\delta_k$ we have to take care of small divisors.
First of all, for each~$k\not\in K_2$ with~$|k|\ge 2$ we can associate to~$\delta_k$ a decomposition of the form
$$\delta_k= \eps_{l_0}^{-1}\eps_{l_1}^{-1}\cdots\eps_{l_q}^{-1},\tag{eqdelta}$$
where~$l_0=k$,~$|k|>|l_1|\ge\cdots\ge|l_q|\ge2$ and~$l_j\not\in K_2$ for all~$j=1, \dots, q$ and~$q\ge 1$. 
If~$k\in K_1$ it is obvious by the definition of~$\delta_k$. If~$k\not\in K_1\cup K_2$, choose a decomposition~$k=k_1+\cdots+k_\nu$ such that the maximum in the expression of~$\delta_k$ is achieved. Obviously,~$k_j$ doesn't belong to~$K_2$ for all~$j=1,\dots, \nu$. We can then express~$\delta_k$ in terms of~$\eps_{k_j}^{-1}$ and~$\delta_{k'_j}$  with~$|k'_j|<|k_j|$. Carrying on this process, we eventually arrive at a decomposition of the form \rf{eqdelta}. Furthermore,
$$\eps_k = |\tilde \lambda^k - \tilde \lambda_{i_k}|, \quad |k|\ge 2,\, k\not\in K_2,$$ 
the index~$i_k$ being chosen in some definite way (of course, if~$k\in K_1$ then~$i_k\in \{s+1, \dots, n\}$).  

\sm The rest of the proof follows closely [P\"o]. For the benefit of the reader, we report here the main steps.

We can define,
$$N^j_m(k), \quad m\ge 2, \quad j\in\{1, \dots, n\},$$
to be the number of factors~$\eps_{l}^{-1}$ in~$\delta_k$, ($l = l_0, \dots, l_q$) satisfying
$$\eps_{l}<\theta\,\tilde\omega(m),~~\hbox{and}~~i_l = j,$$ 
where~$\tilde \omega(m)$ is defined in Definition \rf{De1.0}, and in this notation can be expressed as
$$\tilde\omega(m)= \min_{2\le|k|\le m \atop k\not\in K_2}\eps_k, \quad m\ge2,$$
and~$\theta$ is the positive real number satisfying
$$4\theta=\min_{1\le h\le n}|\tilde\lambda_h| \le 1.$$
The last inequality can always be satisfied by replacing~$f$ by~$f^{-1}$ if necessary. Then we also have~$\tilde\omega(m)\le 2$.

Notice that~$\tilde\omega(m)$ is non-increasing with respect to~$m$ and under our assumptions~$\tilde\omega(m)$ tends to zero as~$m$ goes to infinity. Following [Br], we have the key estimate. 

\sm\thm{Le1.1}{Lemma} {\sl For~$m\ge2$,~$1\le j\le n$ and~$k\not\in K_2$,
$$N^j_m(k)\le \cases{0, &$|k|\le m$,\cr\noalign{\sm}
					\displaystyle {2|k|\over m}-1,&$|k|> m$.}$$}
\sm\proof The proof is done by induction. Since we fix~$m$ and~$j$ throughout the proof, we write~$N$ instead of~$N^j_m$.

For~$|k|\le m$, 
$$\eps_k\ge\tilde\omega(|k|)\ge \tilde\omega(m) > \theta\, \tilde\omega(m),$$
hence~$N(k)=0$.

Assume now that~$|k|>m$. Then~$2|k|/m -1 \ge 1$. If~$k\in K_1$ then, by definition,~$\delta_k = \eps_k^{-1}$, so~$N(k)$ can only be equal to~$0$ or~$1$ and we are done.

Let us suppose~$k\not \in K_1\cup K_2$. Write
$$\delta_k= \eps_k^{-1}
\delta_{k_1}\cdots \delta_{k_\nu}, \quad k=k_1 + \cdots + k_\nu, \quad \nu\ge2,$$
with~$|k|>|k_1|\ge \cdots\ge|k_\nu|$, and consider the following different cases. Observe that~$k - k_1\not\in K_2$, otherwise the other~$k_h$'s would be in~$K_2$. 

\sm{\it Case 1:}~$\eps_k \ge \theta\,\tilde\omega(m)$ and~$i_k$ arbitrary, or~$\eps_k < \theta\,\tilde\omega(m)$ and~$i_k\ne j$. Then
$$N(k) = N(k_1) + \cdots + N(k_\nu),$$
and applying the induction hypotheses to each term we get~$N(k)\le (2|k|/m) - 1$.

\sm{\it Case 2:}~$\eps_k<\theta\,\tilde\omega(m)$ and~$i_k=j$. Then
$$N(k) = 1 + N(k_1) + \cdots + N(k_\nu),$$
and there are three different cases.

{\it Case 2.1:}~$|k_1|\le m$. Then
$$N(k) = 1 < {2|k|\over m} -1,$$
as we want.

{\it Case 2.2:}~$|k_1|\ge|k_2|>m$. Then there is~$\nu'$ such that~$2\le\nu'\le\nu$ and~$|k_{\nu'}|> m\ge |k_{\nu'+1}|$, and we have
$$N(k) = 1 + N(k_1) + \cdots + N(k_{\nu'})\le 1 + {2|k|\over m} - \nu' \le  {2|k|\over m} -1.$$

{\it Case 2.3:}~$|k_1|>m\ge |k_2|$. Then
$$N(k)= 1 + N(k_1),$$
and there are three different cases.

{\it Case 2.3.1:}~$i_{k_1}\ne j$. Then~$N(k_1) = 0$ and we are done. 

{\it Case 2.3.2:}~$|k_1|\le|k|-m$ and~$i_{k_1} = j$. Then
$$N(k) \le 1 + 2\,{|k|-m\over m} -1 < {2|k|\over m} -1.$$

{\it Case 2.3.3:}~$|k_1|>|k|-m$ and~$i_{k_1} = j$. The crucial remark is that~$\eps_{k_1}^{-1}$ gives no contribute to~$N(k_1)$, as shown in the next lemma.

\sm\thm{Le1.2}{Lemma} {\sl If~$k>k_1$ with respect to the lexicographic order,~$k$,~$k_1$ and~$k-k_1$ are not in~$K_2$,~$i_k= i_{k_1} =j$ and
$$\eps_k < \theta\, \tilde\omega(m) \quad\hbox{and}\quad \eps_{k_1}<\theta\, \tilde\omega(m),$$
then~$|k-k_1| = |k| - |k_1| \ge m$.}

\sm\proof Before we proceed with the proof, notice that the equality~$|k-k_1| = |k| - |k_1|$ it is obvious since~$k>k_1$.

Since we are supposing~$\eps_{k_1} = |\tilde\lambda^{k_1}-\tilde\lambda_j| <\theta\,\tilde \omega(m)$, we have
$$\eqalign{|\tilde\lambda^{k_1}| &>|\tilde\lambda_j|-\theta\,\tilde\omega(m)\cr
								&\ge 4\theta - 2\theta = 2\theta.
}$$
Let us suppose by contradiction~$|k-k_1| = |k| - |k_1| < m$. By assumption, it follows that
$$\eqalign{2\theta\,\tilde\omega(m) &> \eps_k+ \eps_{k_1}\cr
							&= |\tilde\lambda^k - \tilde\lambda_{j}| + |\tilde\lambda^{k_1} - \tilde\lambda_{j}|\cr
							&\ge |\tilde\lambda^k - \tilde\lambda^{k_1}|\cr
							&\ge |\tilde\lambda^{k_1}| \,|\tilde\lambda^{k-k_1} - 1|\cr
							&\ge 2\theta\,\tilde \omega(|k-k_1|+1)\cr
							&\ge 2\theta\, \tilde\omega(m),
}$$
which is impossible. \qed

\sm Using Lemma~\rf{Le1.2}, case~$1$ applies to~$\delta_{k_1}$ and we have
$$N(k) = 1 + N(k_{1_1}) + \cdots + N(k_{1_{\nu_{1}}}),$$
where~$|k|> |k_1|> |k_{1_1}| \ge \cdots \ge |k_{1_{\nu_{1}}}|$ and~$k_1 =k_{1_1} + \cdots +k_{1_{\nu_{1}}}$. We can do the analysis of case~$2$ again for this decomposition, and we finish unless we run into case~$2.3.2$ again. However, this loop cannot happen more than~$m+1$ times and we have to finally run into a different case. This completes the induction and the proof of Lemma~\rf{Le1.1}. \qed

\sm Since the origin is a quasi-Brjuno fixed point of order~$s$, there exists a strictly increasing sequence~$\{q_\nu\}_{\nu\ge 0}$ of integers with~$q_0=1$ and such that
$$\sum_{\nu\ge 0} q_\nu^{-1}\log\tilde\omega(q_{\nu +1})^{-1} <\io.\tag{eq8}$$
Since~$\delta_k=0$ for~$k\in K_2$, we have to estimate only
$${1\over |k|}\log\delta_k = \sum_{j=0}^q {1\over |k|} \log \eps_{l_j}^{-1}, \quad k\not\in K_2.$$
By Lemma~\rf{Le1.1}, 
$$\eqalign{{\rm card}\left\{0\le j\le q : \theta\, \tilde\omega(q_{\nu +1}) \le \eps_{l_j} <\theta\, \tilde\omega(q_\nu)\right\} &\le N_{q_\nu}^1(k) + \cdots N_{q_\nu}^n(k) \cr
							&\le {2n|k|\over q_\nu}} $$
for~$\nu\ge 1$. It is also easy to see from the definition of~$\delta_k$ that the number of factors~$\eps_{l_j}^{-1}$ is bounded by~$2|k| - 1$. In particular,
$${\rm card}\left\{0\le j\le q : \theta\, \tilde\omega(q_{1}) \le \eps_{l_j} \right\} \le 2n|k| = {2n|k|\over q_0}. $$
Then,
$$\eqalign{{1\over |k|} \log \delta_k &\le 2n \sum_{\nu\ge 0} q_\nu^{-1} \log(\theta^{-1}\,\tilde\omega(q_{\nu +1})^{-1}) \cr
	&= 2n\left( \sum_{\nu \ge 0} q_\nu^{-1}\log\tilde\omega(q_{\nu +1})^{-1} + \log(\theta^{-1}) \sum_{\nu \ge 0} q_\nu^{-1}\right).}\tag{eq9}$$
Since~$\tilde\omega(m)$ tends to zero monotonically as~$m$ goes to infinity, we can choose some~$\overline{m}$ such that~$1>\tilde\omega(m)$ for all~$m>\overline{m}$, and we get
$$\sum_{\nu\ge\nu_0} q_\nu^{-1} \le {1\over \log \tilde\omega(\overline{m})^{-1}} \sum_{\nu\ge \nu_0} q_\nu^{-1} \log\tilde\omega(q_{\nu +1})^{-1},$$
where~$\nu_0$ verifies the inequalities~$q_{\nu_0 -1}\le \overline{m} < q_{\nu_0}$. Thus both series in parentheses in \rf{eq9} converge thanks to \rf{eq8}. Therefore
$$\sup_k {1\over |k|}\log \delta_k <\io$$ 
and this concludes the proof. \qed

\sm\thm{Rmk2.1}{Remark} Notice that the osculating hypothesis on the~$f$-invariant manifold is necessary. Let us take a look at the following example in~$\C^2$. Let~$f$ be given by
$$\eqalign{&x'=\lambda (1+y)x + x^2\cr
  		   &y'= y }$$
with~$(\lambda,1)$ satisfying the Brjuno condition of order~$1$ (in particular~$\lambda$ is not a root of unity). This germ is not linearizable. In fact, let~$g_y(x)= \lambda(1+y)x + x^2$, so we can write~$f(x,y) = (g_y(x), y)$. A linearization for~$f$ is a germ of holomorphic diffeomorphism~$\psi=(\psi_1, \psi_2)$ fixing the origin, tangent to the identity, and such that
$$\big(g_{\psi_2(x,y)}(\psi_1(x,y)), \psi_2(x,y)\big) = \big(\psi_1(\lambda x, y), \psi_2(\lambda x, y)\big).$$
This last equality implies~$\psi_2\equiv \psi_2(y)$ and~$g_{\psi_2(y)}(\psi_1(x, y)) = \psi_1(\lambda x, y)$. Composing on the right with~$\psi_2^{-1}$ and setting~$h_y(x) = \psi_1(x, \psi_2^{-1}(y))$, we have
$$g_y\big(h_y(x)\big) = h_y(\lambda x).\tag{eqrmk}$$
From \rf{eqrmk} we deduce that~$h_y(0) \in \fix(g_y) = \{0, 1 - \lambda(1+y)\}$. Now,~$h_0(0)= 0$; hence, by continuity~$h_y(0) = 0$ for~$|y|$ small enough, and so~$g_y'(0)h_y'(0) = \lambda h_y'(0)$ for~$|y|$ small enough. But~$h_0'(0)=1\ne 0$; therefore~$\lambda(1+y) = g_y'(0) = \lambda$ for~$|y|$ small enough, which is impossible.
Since~$f$ is not linearizable it cannot admit an osculating invariant manifold of codimension~$1$, even if, obviously, the manifold~$\{x=0\}$ is~$f$-invariant, and~$f$ is linear there. 

\sm\thm{Rmk2.2}{Remark} The reduced Brjuno condition and the hypothesis~$f$ holomorphically linearizable on the osculating manifold are necessary. Consider the following example in~$\C^n$ for~$n\ge 2$. Let~$f$ be a holomorphic diffeomorphism of~$\C^n$, fixing the origin, given by
$$\eqalign{&x_i'=\lambda_i x_i + f_i^i(x,y)\quad\hbox{for}~i=1,\dots, n-1,\cr
  		&y'=\mu y + y^2,}\tag{ex2}$$		
with~${\rm ord}_x(f_i^1)\ge 2$ for every~$i=1, \dots, n-1$,~$(\lambda_1,\dots, \lambda_{n-1}, \mu)$ non resonant, and~$\mu=e^{2\pi \theta}$ with~$\theta\in \R\setminus\Q$ not a Brjuno number. Then~$M=\{x=0\}$ is an osculating manifold of codimension~$n-1$, but~$(\lambda_1,\dots, \lambda_{n-1}, \mu)$  does not satisfy the reduced Brjuno condition (which, since we have no resonances, coincides with the usual Brjuno condition). Furthermore, thanks to Yoccoz's Theorem [Y],~$f|_M$ is not holomorphically linearizable. This germ is not holomorphically linearizable. 
In fact, assume by contradiction that~$\psi$ is a holomorphic linearization. Then~$\tilde M=\psi(M)=\{\psi^{-1}_1(\tilde x, \tilde y)=0, \dots, \psi^{-1}_{n-1}(\tilde x,\tilde y)=0\}$ is an osculating manifold of codimension~$n-1$ for~$\tilde f(\tilde x,\tilde y)= \psi\circ f\circ \psi^{-1} \equiv \diag(\lambda_1,\dots, \lambda_{n-1},\mu)(\tilde x,\tilde y)$. Thanks to the implicit function Theorem there exist~$n-1$ holomorphic functions~$\chi_1(\tilde y),\dots, \chi_{n-1}(\tilde y)$, such that~$\tilde M= \{\tilde x_1= \chi_1(\tilde y),\dots, \tilde x_{n-1}= \chi_{n-1}(\tilde y)\}$. The~$\tilde f$-invariance of~$\tilde M$ yields
$$\lambda_i \chi_i(\tilde y) = \chi(\mu\tilde y) \quad \hbox{for}~i=1, \dots, n-1,$$
and this is equivalent, writing~$\chi_i(\tilde y)= \sum_{m\ge 1} \chi^i_m \tilde y^m$, to
$$\sum_{m\ge 1}\lambda_i \chi^i_m \tilde y^m = \sum_{m\ge 1} \chi^i_m \mu^m\tilde y^m,$$
which implies~$\chi^i_m\equiv 0$ for every~$i=1, \dots, n-1$ and~$m\ge 0$, because~$(\lambda_1, \dots, \lambda_{n-1}, \mu)$ is not resonant. Then~$\tilde M= \{\tilde x = 0\}$ and, since~$\tilde f|_{\tilde M}$ is linear, we have a holomorphic linearization of~$f|_M$, contradiction.

\sect Final remarks

We can obtain many of the result recalled in the Introduction as corollaries of our Theorems. If there are no resonances Theorem \rf{Te1.0} with~$s=n$ yields Theorem \rf{Te0.1}. If there are no resonances and the origin is an attracting [resp., repelling] fixed point then Theorem \rf{Teorema} with~$s=n$ yields Theorem \rf{Te0.3} because the Brjuno condition is automatically satisfied. 

Our result can be also compared with the following result obtained by Nishimura in [N] (the statement is slightly different from the original one presented in [N] but perfectly equivalent): 

\sm\thm{TeNishimura}{Theorem}(Nishimura, 1983 [N]) {\sl Let~$f$ be a germ of holomorphic diffeomorphism of~$\C^n$, fixing the origin~$O$. Assume that~$Y$ is a complex manifold through~$O$ of codimension~$s$ pointwise fixed by~$f$. In coordinates~$z=(x,y)$ in which~$Y=\{x=0\}$ we can write~$f$ in the form
$$\casi{x_i'= \sum_{k=1}^s C_{ik}(y) x_k  + f^1_i(x,y) &\hbox{for}~$i=1,\dots, s$,\cr\noalign{\sm}
 			\displaystyle{y_j'= y_j + f^2_j(x,y)} &\hbox{for}~$j=1,\dots, r$, }$$
with~${\rm ord}_x(f_i^1)\ge 2$ and~${\rm ord}_x(f_j^2)\ge 1$. If for each point~$p\in Y$ the eigenvalues~$\{\lambda_1(p), \dots, \lambda_s(p)\}$ of the matrix~$C(p) = \big(C_{jk}(p)\big)$ have modulus less than~$1$ and have no resonances, then there exists a unique holomorphic change of coordinates~$\psi$, defined in a neighbourhood of~$Y$, tangent to the identity such that
$$f\circ \psi = \psi\circ L,$$
where~$L$ is the germ
$$\casi{ x_i'= \sum_{k=1}^s C_{ik}(y) x_k &\hbox{for}~$i=1,\dots, s$,\cr\noalign{\sm}
  			\displaystyle{y_j'=y_j} &\hbox{for}~$j=1,\dots, r$.}$$
}

\sm The hypotheses of Nishimura are slightly different from ours, and, in fact,  he does not prove a true linearization theorem. However, his result becomes a linearization result when~$C(y)$ is a constant matrix, which is equivalent to requiring that~$Y$ is an osculating fixed manifold. In this situation our result can be seen as a generalization of Theorem \rf{TeNishimura} in the case of~$\d f_O$ diagonalizable. In fact while he needs an osculating fixed manifold and a strong hypothesis on the modulus of the eigenvalues, we only need an osculating manifold on which our germ is holomorphically linearizable and the origin as a quasi-Brjuno fixed point of order~$s$.

\me Recently, Rong [R] proved the following result

\sm\thm{TeRong}{Theorem}(Rong, 2006 [R]) {\sl Let~$f$ be a germ of holomorphic diffeomorphism of~$\C^n$, fixing the origin with~$\d f_O = \diag(\Lambda_s, I_r)$, where~$\Lambda_s= \diag(\lambda_1, \dots, \lambda_s)$ with~$\lambda_j = e^{2\pi i \theta_j}$,~$\theta_j\in\R\setminus\Q$. Let~$M$ be a pointwise fixed complex manifold through~$O$ of codimension~$s$. Choose local coordinates~$(x,y)$ centered in~$O$ such that~$M=\{x=0\}$. For any~$p\in M$, write~$\d f_p=\left(\matrix{\Lambda_s(y)& O\cr\star  &I_r }\right)$. Assume that~$\Lambda_s(y)\equiv\Lambda_s$ for all~$p\in M$. If the~$\lambda_j$'s satisfy the Brjuno condition, then there exists a local holomorphic change of coordinates~$\psi$ such that
$$ f\circ \psi = \psi\circ \Lambda,$$
where~$\Lambda$ is the linear part of~$f$.}

\me This result too can be seen as a particular case of Theorem \rf{Teorema}. In fact, if we are in the hypotheses of Rong, our hypotheses are automatically verified:~$M$ is an osculating fixed manifold thanks to the hypothesis~$\Lambda_s(y)\equiv\Lambda_s$ for all~$p\in M$, and the hypotheses on the eigenvalues follow immediately from Remarks \rf{Re0.2} and \rf{Re1.1}. 

\sm What it is new in our result is that we are not assuming anything on the modulus of the eigenvalues, so we are really dealing with the mixed case. In fact we are allowing cases in which there are some eigenvalues with modulus greater than 1, some eigenvalues with modulus~$1$, and the remaining eigenvalues with modulus less than~$1$. Finally, our Theorem applies in cases not covered by the previous results, as shown by Remark \rf{Re1.1}.

%

\vbox{\vskip.75truecm\advance\hsize by 1mm
	\hbox{\centerline{\sezfont References}}
	\vskip.25truecm}\nobreak
\parindent=40pt

\bib{A} {\sc Abate, M.:} {\sl Discrete local holomorphic dynamics,} in ``Proceedings of 13th Seminar of Analysis and its Applications, Isfahan, 2003'', Eds. S. Azam et al., University of Isfahan, Iran, 2005, pp. 1--32.

\bib{Ar} {\sc V.I. Arnold:} ``Geometrical methods in the theory of ordinary differen\-tial equations'', Springer-Verlag, Berlin, 1988.

\bib{BCS} {\sc Bao, D., Chern, S.-S., Shen, Z.:} ``An introduction to Riemann-Finsler geometry'', Graduate Texts in Mathematics {\bf 200}, Springer-Verlag, New York, 2000.

\bib{BER} {\sc Baouendi, M. S., Ebenfelt, P., Rothschild, L. P.:} {\sl Dynamics of the Segre varieties of a real submanifold in complex space,} J. Algebraic Geom. {\bf 12} (2003), no. 1, pp. 81--106.

\bib{BMR} {\sc Baouendi, M. S., Mir, N., Rothschild, L. P.:} {\sl Reflection ideals and mappings between generic submanifolds in complex space,} J. Geom. Anal. {\bf 12} (2002), no. 4, pp. 543--580.

\bib{B} {\sc Bracci, F.:} {\sl Local dynamics of holomorphic diffeomorphisms,} Boll. UMI (8), 7--B (2004), pp. 609--636.

\bib{Br} {\sc Brjuno, A. D.:} {\sl Analytic form of differential equations,} Trans. Moscow Math. Soc., {\bf 25} (1971), pp. 131--288; {\bf 26} (1972), pp. 199--239.

\bib{C} {\sc Carletti, T.:} {\sl Exponentially long time stability for non-linearizable analytic germs of~$(\C^n, 0)$,} Annales de l'Institut Fourier (Grenoble), {\bf 54} (2004), no. 4, pp. 989--1004.

\bib{D} {\sc Dulac, H.:} {\sl Recherches sur les points singuliers des \'equationes diff\'erentielles}, J. \'Ecole polytechnique II s\'erie cahier IX, (1904), pp. 1--125.

\bib{G} {\sc Gray, A.:} {\sl A fixed point theorem for small divisors problems,}  J. Diff. Eq., {\bf 18} (1975), pp. 346--365.

\bib{M} {\sc Marmi, S.:} ``An introduction to small divisors problems'', I.E.P.I., Pisa, 2003.

\bib{N} {\sc Nishimura, Y.:} {\sl Automorphismes analytiques admettant des sous-vari\'et\'es de points fix\'es attractives dans la direction transversale,} J. Math. Kyoto Univ., {\bf 23--2} (1983), pp. 289--299.

\bib{P} {\sc Poincar\'e, H.:} ``\OE uvres, Tome I'', Gauthier-Villars, Paris, 1928, pp. XXXVI--CXXIX.

\bib{P\"o} {\sc P\"oschel, J.:} {\sl On invariant manifolds of complex analytic mappings near fixed points,} Exp. Math., {\bf 4} (1986), pp. 97--109.

\bib{R} {\sc Rong, F.:} {\sl Linearization of holomorphic germs with quasi-parabolic fixed points,} to appear in Ergodic Theory Dynam. Systems. 

\bib{R\"u} {\sc R\"ussmann, H.:} {\sl Stability of elliptic fixed points of analytic area-preserving mappiongs under the Brjuno condition,} Ergodic Theory Dynam. Systems, {\bf 22} (2002), pp. 1551--1573.

\bib{S} {\sc Sternberg, S.:} {\sl Infinite Lie groups and the formal aspects of dynamical systems,} J. Math. Mech., {\bf 10} (1961), pp. 451--474.

\bib{Y} {\sc Yoccoz, J.-C.:} {\sl Th\'eor\`eme de Siegel, nombres de Bruno et polyn\^omes quadratiques,} Ast\'erisque {\bf 231} (1995), pp. 3--88.

\bye